\documentclass[12pt]{amsart}

\usepackage{amssymb, amsfonts}

\theoremstyle{plain}

\newcommand{\beq}{\begin{equation}}
\newcommand{\eeq}{\end{equation}}

\newcommand{\rti}{\textbf{$r\rightarrow \infty$}}
\newcommand{\tti}{\textbf{$t\rightarrow \infty$}}

\newtheorem{cor}{Corollary}
\newtheorem{theor}{Theorem}
\newtheorem{prop}{Proposition}
\newtheorem{rem}{Remark}

\numberwithin{equation}{section}

\begin{document}

\begin{center}
{\large {\bf Valiron-Titchmarsh' and Related Theorems for Subharmonic
Functions in $\mathbb{R}^n$ With Masses on a Half-Line}} \\
\vspace{1cm}

\vspace{2cm}

Alexander I. Kheyfits \\

\medskip

Bronx Community College and the Graduate Center \vspace{.1cm} \\
of the City University of New York, USA \vspace{.1cm} \\

\medskip

\email{akheyfits@gc.cuny.edu}

\vspace{1cm}

\end{center}

\vspace{1cm}

\noindent {\bf {Abstract}.} The Valiron-Titchmarsh theorem on asymptotic behavior of entire functions with negative zeros is extended to subharmonic functions with the Riesz masses on a ray in $\mathbb{R}^n,\; n\geq 3$. We also show that the existence of the limit $\lim_{\rti} \frac{\log u(r)}{N(r)}$, where $N(r)$ is the averaged counting function of a subharmonic function $u$ with the associated masses on a ray implies the regular asymptotic behavior separately for $u$ and for $N$.
\vspace{3cm}

\begin{center}
\today
\end{center}

\vspace{1cm}

\footnoterule \vspace{.5cm}



\newpage

\section{Introduction and Statement of Results}

In 1913, Georges Valiron published his well-known memoir in \emph{Annales de la facult\'{e} des sciences de Toulouse} \cite{Val}, which became his dissertation the next year. A century later, this work is still regularly cited. The well-known result of \cite{Val} is the sequel Tauberian theorem about entire functions with roots on a half-line. An elegant succinct account of the current research about this class of functions has been recently given by Drasin \cite{Dra2}.  \\

\noindent\textbf{Theorem} (Valiron\footnote{In fact, Valiron proved the theorem for the entire functions of proximate order and so-called directed products.} \cite[p. 237]{Val}). \emph{Let} $f(z), z=r e^{\imath \theta}$, \emph{be an entire function of non-integer order} $\rho$ \emph{and finite type with negative zeros, and} $n(t)$ \emph{be the counting function of its zeros, that is, the number of the zeros of} $f$ \emph{in the closed disk} $\{|z| \leq t\}$.

\emph{If there exists the limit over the positive ray}
\[\lim_{\rti} r^{-\rho} \log f(r) = \frac{\pi \mathbf{\Delta}}{\sin \pi \rho}, \]
\emph{where} $\mathbf{\Delta}$ \emph{is a constant, then there exists the limit } 
\beq \lim_{\tti} t^{-\rho} n(t) = \mathbf{\Delta}. \eeq

The positivity of the counting function is the Tauberian condition of the theorem.  \\

A proof of the converse Abelian assertion is straightforward.  \\

\emph{If there exists the limit (1.1), then there exist the limits}
\beq \lim_{\rti} r^{-\rho} \log f(r e^{\imath \theta}) = \frac{\pi \mathbf{\Delta}}{\sin \pi \rho}e^{\imath \rho \theta}   \eeq
\emph{for any} $\theta \in (-\pi, \pi)$, \emph{uniformly in any sector} $-\pi+\delta \leq \theta \leq \pi -\delta$, $0<\delta <\pi$, \emph{thus} $f$ \emph{is a function of the completely regular growth in the plane}.  \\

These results translate assumption (1.1) on the regular asymptotic behavior of the counting function into conclusion (1.2) about the regular growth of the entire function, and vice versa. If some information about the \emph{joint} regular behavior of an entire function and its counting function is known, for example, if we know the existence of the limit $\lim_{\rti} \frac{\log f(r)}{n(r)}$, then under some conditions one can derive certain conclusions on the \emph{separate} asymptotic behavior of the function and the counting function of its zeros \cite{Dra1}. We give a result of this nature for subharmonic functions with the associated masses on a ray as well.   \\

The Valiron theorem was later independently proved by Titchmarsh \cite{Tit}, thus it is often referred to as the Valiron-Titchmarsh theorem on entire functions with negative zeros; its current exposition can be found, for example, in \cite[Lect. 12-13]{Lev}. Delange proved that the positive $x-$axis in (1.2) can be replaced by any ray
\beq \phi \neq \frac{2k+1}{2\rho}\pi \mbox{ with an integer } k,\;-\rho -1/2<k<\rho -1/2, \eeq
see \cite[p. 200, Theor. 5.8.1.2]{Azb}. The result was extended to holomorphic and subharmonic functions in $\mathbb{R}^2$ and in plane angles, with associated masses on several rays or on logarithmic spirales, and to functions with many-term asymptotics, see \cite{AgrLo, Az97, Azb, Khe1, Khe2} and the references therein. Agranovich \cite{Agr} considered the Abelian part of the theorem for subharmonic functions with many-term asymptotics, with masses on one ray in $\mathbb{R}^n,\; n\geq 3$.

In this paper we extend the Valiron-Titchmarsh theorem to subharmonic functions with masses on a ray in $\mathbb{R}^n,\; n\geq 3$. We use the approach based on the General Tauberian Theorem by Wiener \cite[Chap. V]{PaW}. It invokes an integral representation of the indicators of the subharmonic functions under consideration through the associated Legendre functions $P_{\nu}^{\mu}(\cos \psi)$ of the first kind. To deduce this, likely new representation, we compute in the Appendix the Mellin integral transform of the Weierstrass primary kernel. In this section we fix notation, following \cite{HaKe} and \cite{GO}, and state our results. The proofs are given in Section 2.   \\

Introduce in $\mathbb{R}^n=\{x=(x_1,\ldots ,x_n)\}$ spherical coordinates $x =(r,\theta),\;r=|x|$, $\theta =(\theta_1, \ldots, \theta_{n-1})$, such that $x_1=r\cos \theta_1$. Here $0\leq \theta_1 \leq \pi$ and $0\leq \theta_k \leq 2\pi$ for $k=2,3,\ldots ,n-1$. The ball of radius $t$ centered at the origin of $\mathbb{R}^n$ is denoted as $B_t$, $\overline{B_t}$ being its closure.

Let $u(x)$ be a subharmonic function of finite non-integer order $\rho$ in $\mathbb{R}^n,\; n\geq 3$, such that its Riesz associated measure $\mu$ is supported by the negative $x_1-$axis, and let
\[n(t) = t^{2-n}\mu(\overline{B_t}) \]
be the counting function of the measure $\mu$. By the Hadamard representation theorem \cite[Sect. 4.2.2]{HaKe}
\beq u(x)= \int_{\mathbb{R}^n} K_q(x, y)d\mu(y) +u_0(x). \eeq
Here $K_q$ is the Weierstrass canonical kernel,
\[K_q(x, y) = -(r^2+t^2-2tr \cos \psi)^{-\frac{n-2}{2}}
+ t^{2-n}\sum_{j=0}^q \left(\frac{r}{t}\right)^j G^{\frac{n-2}{2}}_j (\cos \psi), \]
the integral in (1.4) is convergent uniformly on any compact set and absolutely at the vicinity of infinity, $u_0$ is a harmonic polynomial of degree at most $q$, $q=E(\rho)$ being the integer part of $\rho,\; 0\leq q<\rho <q+1$. Also, let $t=|y|$, $\psi=\widehat{(x,y)}$ be the angle between vectors $x$ and $y$, and $G^{(n-2)/2}_j (\cos \psi)$ the Gegenbauer polynomials given by the generating function
\[(1-2t \cos \psi + t^2)^{-\frac{n-2}{2}} = \sum_{j=0}^{\infty} G^{\frac{n-2}{2}}_j (\cos \psi)\; t^j.\]
Since we are interested in the asymptotic properties of $u$, which are not affected by the polynomial $u_0$, we assume hereafter
that $u_0=0$. Due to the same reason, we suppose, without loss of generality, that the closed unit ball $\overline{B_1}$ is free of masses of $u$, thus $n(t) = 0$ for $0\leq t \leq 1$.

For a vector $y=(t,\theta)$ with $\theta_1=\pi$, all the other angular coordinates are undefined, and we represent such a vector as $y=(t,\pi)$. If a measure $\mu$ is distributed over the negative $x_1-$axis, the angle between the vectors $x=(r,\theta_1,\ldots , \theta_{n-1})$ with the latitude $\theta_1$, and $y =(t,\pi)$ in representation (1.4) is $\psi=\pi -\theta_1$ for any $\theta_2,\ldots ,\theta_{n-1}$, and integral (1.4) becomes
\beq u(x)= \int_0^{\infty} K_q(x, (t,\pi))d\left(t^{n-2}n(t)\right).  \eeq
Moreover, due to our assumptions, integral (1.5) can be written, see \cite[Eq-ns (9) and (11)]{GO}, as
\beq u(x)= \int_0^{\infty} t^{2-n} h_n\left(\frac{r}{t},\theta_1 ,q \right) d\left(t^{n-2}n(t)\right),  \eeq
where the kernel is
\[h_n(s,\theta_1,q)=-\left(1+2s \cos \theta_1+s^2\right)^{\frac{2-n}{2}}+\sum^q_{j=0}(-s)^j G_j^{\frac{n-2}{2}}(\cos \theta_1).\]
\begin{rem} It is worth mentioning that in the case under consideration the function $h_n(s,\theta_1,q)$ and the right-hand side of formula (1.5) depend only on the angle $\theta_1$, but not upon the other angular coordinates $\theta_2,\ldots ,\theta_{n-1}$ of the point $x$.   \end{rem}

We first state the Abelian result regarding the subharmonic functions given by (1.5)-(1.6).

\begin{prop} Let $u$ be a subharmonic function in $\mathbb{R}^n,\; n\geq 3$, of non-integer order $\rho$ and finite type, whose Riesz masses are distributed over the negative $x_1-$axis. If there exists the limit
\beq \lim_{\tti} t^{-\rho} n(t) = \mathbf{\Delta}, \eeq
then for any $x=(r,\theta)$, with $0\leq \theta_1< \pi $, there exists the limit
\beq \lim r^{-\rho}u(x) = H(\theta_1). \eeq
This limit, that is, the indicator function of the subharmonic function $u$, is given by
\beq H(\theta_1) =(\rho +n-2) \mathbf{\Delta} \int_0^{\infty} s^{-\rho -1} h_n(s,\theta_1,q) ds.  \eeq

Moreover, the indicator can be expressed through the associated Legendre spherical functions of the first kind 
\[P_{\nu}^{\mu}(\cos \theta_1)\] 
on the cut \cite[Chap. 3]{BE}, as
\[H(\theta_1)=\frac{\pi 2^{\frac{n-3}{2}}\Gamma(\frac{n-1}{2})\prod_{k=1}^{n-2}(\rho+k)\; \mathbf{\Delta}}{(n-3)! \sin (\pi \rho)
\left(\sin \theta_1\right)^{\frac{n-3}{2}}}\; P^{\frac{3-n}{2}}_{-\rho-\frac{n-1}{2}}(\cos \theta_1). \]
Using the known property of the Legendre functions, $P^{\mu}_{\nu}=P^{\mu}_{-\nu-1}$ \cite[Sect. 3.3.1(1)]{BE}, the latter can be rewritten as
\beq H(\theta_1)=\frac{\pi 2^{\frac{n-3}{2}} \Gamma(\frac{n-1}{2}) \prod_{k=1}^{n-2}(\rho+k)\; \mathbf{\Delta}}{(n-3)! \sin (\pi \rho )
\left(\sin \theta_1\right)^{\frac{n-3}{2}}}\; P^{\frac{3-n}{2}}_{\rho +\frac{n-3}{2}}(\cos \theta_1). \eeq
Equation (1.10) holds good for $\theta_1=\pi$ as well, since in this case both its sides are equal to $-\infty$.
\end{prop}

The fact that (1.7) implies (1.8) is not new, it was established by Azarin \cite{Az} as early as in 1961, together with an integral representation of the indicator of subharmonic functions in $\mathbb{R}^n,\; n\geq 3$, in terms of the fundamental solutions of the Legendre differential equation
\[t(1-t)y''(t)+(n-1)\left(\frac{1}{2}-t\right)y'(t)+\rho (\rho+n-2)y(t)=0.\]
The occurrence of the associated Legendre functions in problems like ours was mentioned in \cite[p. 160]{HaKe}, however, without explicit formulas.

A new feature of our result is the explicit representation (1.10) of the indicator in terms of the associated Legendre functions of the first kind, leading to a precise description of the zero sets of the indicators of subharmonic functions at question. We need this description to apply Wiener's tauberian theorem.

The zeros of the spherical functions have been carefully studied and (before the computers became omnipresent) tabulated. It is known in particular, that if $\nu$ is not real, then $P^{\mu}_{\nu}(\cos \beta)$ is never zero \cite[p. 403]{Hob}, while for any real $\mu$ and $\nu$ the equation $P^{\mu}_{\nu}(\cos \beta)=0$ has only finitely many roots. In our case $\mu=(3-n)/2$ and $\nu=\rho +(n-3)/2$, hence the function
\[f(\beta)=P^{(3-n)/2}_{\rho +(n-3)/2}(\cos \beta) \]
has finitely many real zeros on $(0, \pi)$.

More precisely \cite[p. 386-388]{Hob}, the equation $f(\beta)=0$ has $E(\rho+1)=q+1$ zeros in the interval $(0,\pi)$. In particular, if $0<\rho<1$, then for any dimension $n$ there is only one root, which was observed in \cite[p.161]{HaKe}. For instance \cite{Ba}, if $n=3$ and $\rho \approx 0.5$, then the only root $\beta_1 \approx 130^{\circ}$; if $n=5$ and $\rho \approx 0.5$, then the unique root $\beta_1 \approx 115^{\circ}$.

Let $\Theta_n(\rho)=\{\beta^n_1,\ldots , \beta^n_{q+1}\}$ stand for the set of the roots of the equation $f(\beta)=0$. Similarly to condition (1.3), our results include the restriction $\phi \not\in \Theta_n(\rho)$. The occurrence of a finitely many exceptional rays in the Valiron-Titchmarsh theorems was discovered by Delange \cite{Del} and studied in detail in \cite{Az97}. \\

Now we state the Tauberian counterpart of Proposition 1.

\begin{theor} Let $u$ be a subharmonic function in $\mathbb{R}^n,\; n\geq 3$, of non-integer order $\rho$ and finite type, whose Riesz masses are distributed over the negative $x_1-$axis. Let $\phi,\; 0\leq \phi<\pi$, be the angle between a vector $x\in \mathbb{R}^n$ and the positive $x_1-$axis. If $\phi \not\in\Theta_n(\rho)$ and the limit
\beq \lim r^{-\rho}u(x) = H(\phi) \eeq
exists, then there exists the limit
\[\lim_{\tti} t^{-\rho} n(t) = M(g;0) H(\phi) \]
where
\[M(g,0) =\frac{2^{(n-3)/2}\Gamma((n-2)/2) \sin (\pi \rho)\left(\sin \phi\right)^{(n-3)/2}}{\pi^{3/2} \prod ^{n-3}_{k=1}(\rho +k) P^{(3-n)/2}_{\rho +(n-3)/2}(\cos \phi)} \]
and therefore, by Proposition 1, there exist the limits for any $\theta_1,\; 0\leq \theta_1<\pi$,
\[\lim r^{-\rho}u(x) = H(\theta) \equiv H(\theta_1),\;\;  x=(r,\theta_1,\ldots, \theta_{n-1}), \]
with
\beq H(\theta_1)=\left(\frac{\sin \phi}{\sin \theta_1}\right)^{\frac{n-3}{2}} \frac{P^{(3-n)/2}_{\rho +(n-3)/2}(\cos \theta_1)} {P^{(3-n)/2}_{\rho +(n-3)/2}(\cos \phi)} H(\phi).   \eeq

If $\phi \in\Theta_n(\rho)$, the conclusion fails as an example below shows.
\end{theor}
\begin{rem} In $n=2$ there exists an entire function, whose indicator vanishes on finitely many rays, and the function has the completely regular growth on these and only these rays \cite[p. 161]{Lev0}.
\end{rem}
\begin{rem} If $n=2$ the Weierstrass primary kernel is different from that in (1.4), therefore in the plane case our proof is invalid. Nonetheless, if $n=2$, we have \cite{BE}
\[P^{1/2}_{-\rho -1/2}(\cos \alpha)=\sqrt{2/(\pi \sin \alpha)} \cos (\rho \alpha),\; 0<\alpha <\pi,\]
thus in this case formula (1.12) becomes the known one, \cite[Theor. 1']{Khe1},
\[H(\theta_1)=\frac{H(\phi)}{\cos \rho \phi}\cos \rho \theta_1.\]  \end{rem}

Finally we consider the asymptotic behavior of the ratios $u(r,\theta)/n(r)$ for a subharmonic function $u$ satisfying the conditions of Theorem 1. The precise bounds for this ratio were obtained by Gol'dberg and Ostrovskii \cite{GO}:  \\

\emph{For a Weierstrass canonical integral} $v$ \emph{of noninteger order} $\rho,\; q=E(\rho)$, \emph{with Riesz masses on the negative} $x_1-$\emph{half-axis}
\[\liminf_{\rti}\frac{v(r,\theta)}{n(r)}\leq (\rho+n-2)\int_0^{\infty}u^{-1-\rho}h_n(u,\theta_1,q)du\leq \limsup_{\rti}\frac{v(r,\theta)}{n(r)},\]
\emph{where} $\theta$ \emph{is any point of the unit sphere with} $0\leq \theta_1 <\pi$.   

Gol'dberg and Ostrovskii showed that both upper and lower bounds in their theorem are exact. Theorem 2 below shows that if we replace the function $r^{\rho}$ with $r^{\rho (r)}$ with a proximate order $\rho(r)$, this gives all the functions the bounds at \cite{GO} are attained at. \\

We are interested in the asymptotic behavior of the ratio $v(r,\theta)/n(r)$. For the entire functions, the following theorem was proved by Drasin \cite{Dra1}. It shows, in particular, that the proximate orders cannot be avoided in this case. \\

\emph{If} $f$ \emph{is an entire function of order} $\rho,\; 0<\rho<1$, \emph{with all zeros real and negative, then either of the conditions}
\[\log M(r)/n(r) \rightarrow L>0, \; \rti ,\]
\emph{or}
\[\log M(r)/N(r) \rightarrow L/\rho >0, \; \rti ,\]
\emph{imply that}
\[\log M(r)=r^{\rho} \psi(r),\]
\emph{where} $\rho$ \emph{is determined by} $L=\pi /\sin (\pi \rho)$ \emph{and} $\psi$ \emph{is a slowly varying function, that is,} $\psi(\sigma r)/\psi(r)\rightarrow 1$ \emph{as} $\rti$ \emph{for each fixed} $\sigma >0$.  \\

Instead of the counting function $n$, here it is more convenient to use the function $\overline{n}(t)=\mu(\overline{B_t})$ and its average 
\[N(r)=(n-2)\int_0^r t^{1-n}\overline{n}(t) dt.\]

We first state the following immediate corollary of Theorem 1.

\begin{cor} Under the conditions of Proposition 1 or Theorem 1, there exist the limits, where $x=(r,\theta_1,\ldots,\theta_{n-1})$,
\[\lim_{\rti}\frac{u(x)}{n(r)}=\frac{\pi 2^{\frac{n-3}{2}} \Gamma(\frac{n-1}{2}) \prod_{k=1}^{n-2}(\rho+k)}{(n-3)!\sin(\pi \rho)
\left(\sin \theta_1\right)^{\frac{n-3}{2}}}\; P^{\frac{3-n}{2}}_{\rho +\frac{n-3}{2}}(\cos \theta_1)  \]
and
\[\lim_{\rti}\frac{u(x)}{N(r)}= \frac{\pi 2^{\frac{n-3}{2}} \Gamma(\frac{n-1}{2}) \prod_{k=0}^{n-2}(\rho+k)}{(n-2)!\sin(\pi \rho)
\left(\sin \theta_1\right)^{\frac{n-3}{2}}}\; P^{\frac{3-n}{2}}_{\rho +\frac{n-3}{2}}(\cos \theta_1). \]   \\
\end{cor}

Our last result is a Tauberian theorem for the subharmonic functions, similar to the cited theorem of Drasin. Let \cite[Eq-n (4.5.16)]{HaKe}
\[P_n(r,t,\theta_1)  \]
\[=rt^{n-2}\left((n-1)r^2\cos \theta_1+rt[n+(n-2)\cos ^2 \theta_1]+(n-1)t^2 \cos \theta_1\right). \]

\begin{theor} Let $u$ be a subharmonic function in $\mathbb{R}^n,\; n\geq 3$, of order $\rho,\;0<\rho <1$, whose Riesz masses are distributed over the negative $x_1-$axis. If the limit
\[\lim_{\rti}\frac{u(r)}{n(r)}=\overline{\mathbf{\Delta}}  \]
exists, then
\[N(r)= r^{\rho} \psi_1(r),  \]
and
\[u(r)=r^{\rho} \psi(r), \]
where $\psi$ and $\psi_1$ are slowly varying functions for $0<r<\infty$ and the order $\rho$ satisfies the equation
\[\overline{\mathbf{\Delta}}=\frac{ \Gamma(n-1-\rho)}{(n-2)! \Gamma(1-\rho)} \frac{\pi \rho}{\sin (\pi \rho)}.\]
Moreover,
\[u(x)=B(\theta) r^{\rho} \psi(r),\]
where
\[B(\theta)=B(\theta_1)=\int_0^{\infty}\frac{t^{\rho} P_n(t,1,\theta_1)\psi_1(t) dt}
{\left(t^2+2t\cos \theta_1+1\right)^{(n+2)/2}}.\]
\end{theor}
\begin{rem} In the assumptions of the theorem, we can replace the positive $x_1-$axis with any ray $0\leq \theta_1 \leq \pi/2$, but the corresponding transcendental equation is cumbersome, and we leave it out.   \end{rem}

\section{Proofs}

\noindent\textbf{Proof of Proposition 1.} By condition (1.1), $n(t) = \Delta t^{\rho} + t^{\rho} \varepsilon(t)$, with $\lim_{\tti}\varepsilon(t)=0$, thus
\[u(x) = u_1(x) + u_2(x), \]
where
\[u_1(x)= \int_0^{\infty} K_q(x, (t,\pi))d\left(t^{n-2} \Delta t^{\rho}\right)  \]
\[=(\rho +n-2) \mathbf{\Delta} \int_0^{\infty} t^{\rho +n-3} K_q(x, (t,\pi))dt  \]
and $u_2=u-u_1$.

Due to (1.6),
\[u_1(x) =(\rho +n-2) \mathbf{\Delta} \int_0^{\infty} t^{\rho -1} h_n \left(\frac{r}{t}, \theta_1, q \right) dt  \vspace{.3cm}\]
\[=(\rho +n-2) \mathbf{\Delta} \left(\int_0^{\infty} s^{-\rho -1} h_n \left(s, \theta_1, q \right) ds \right) r^{\rho},  \]
thus giving (1.9)-(1.10).

Next we estimate the second term $u_2$, which is a subharmonic function of the non-integer order $\rho$. We need the known bound of the kernel $h_n$, see for example, \cite{GO},
\beq \big| h_n(u,\theta,q) \big| \leq C\min \left(u^q; u^{q+1}\right) \eeq
for all $u>0$ with a constant $C$ depending only on $n$ and $q$.

Since the associated measure of $u_2$ has the minimal type with respect to the non-integer $\rho$, the estimates in \cite[Chap. 2]{Ron} or \cite[Chap. 4]{HaKe} imply that the function $u_2$ itself has zero type with respect to $\rho$, thus proving (1.8). \\

To verify (1.10) when $\theta_1=\pi$, we notice that in this case the integral in (1.10) is divergent to $-\infty$. On the other hand, by making use of the known asymptotic formulas for the associated Legendre function \cite[Sect. 3.9.1]{BE}, we straightforwardly find that as $\theta_1 \uparrow \pi$,
\[H(\theta_1)\approx -\frac{(\rho +n-2)\left(\Gamma \left(\frac{n-3}{2}\right)\right)^2 \mathbf{\Delta} }{2(n-4)!} \frac{1}{(\cos (\theta_1 /2))^{n-3}} \searrow -\infty\]
for $n\geq 4$, and
\[H(\theta_1)\approx (\rho+1) \mathbf{\Delta} \left[2\log (\cos (\theta_1 /2)) +\gamma +2\psi(-\rho)-\pi\cot (\pi \rho) \right] \searrow -\infty \]
if $n=3$. Here $\gamma$ is the Euler-Mascheroni constant and $\psi$ the logarithmic derivative of the $\Gamma-$function.  \hfill   $\qed$   \\

Before proving Theorem 1, we formulate the following variant of the General Tauberian Theorem of Wiener \cite[Sect. 12.8, Theor. 233]{Hard}.   \\

\noindent\textbf{Theorem W.} \emph{Let a function} $g\in L(0,\infty)$ be such that
\[\int_0^{\infty}g(t)t^{-\imath x}dt\neq 0 \]
\emph{for all real} $x$. \emph{Let a bounded real function} $f$ \emph{be slowly decaying, that is,}
\[\underline{\lim}_{x\rightarrow \infty}\left(f(y)-f(x)\right) \geq 0 \]
\emph{as} $y>x$ \emph{and} $\frac{y}{x}\rightarrow 1$. \emph{If there exists the limit}
\[\lim_{x\rightarrow \infty}\frac{1}{x}\int_0^{\infty} g \left(\frac{t}{x}\right)f(t)dt=l\int_0^{\infty}g(t)dt, \]
\emph{then there exists the limit
\[\lim_{x\rightarrow \infty}f(x)=l \]
with the same constant} $l$.

\begin{rem} A proof of the Valiron-Titchmarsh theorem for analytic functions in $\mathbb{R}^2$ can be based on Montel's theorem \cite[p. 464-465]{Lev0}, however, Montel's theorem is not valid for subharmonic functions \cite{Bear}.
\end{rem}

\noindent\textbf{Proof of Theorem 1}. We again represent $u(x)$ by (1.6) and derive the equation
\beq\frac{u(r,\theta)}{r^{\rho}}=\frac{1}{r}\int_0^{\infty} \left(\frac{r}{t}\right)^{1-\rho}h_n\left(\frac{r}{t},\theta_1,q\right)d\alpha(t), \eeq
where
\[d\alpha(t)=t^{3-n-\rho}d\left(t^{n-2}n(t)\right).\]
Set in (2.2) $\theta_1=\phi$ from (1.11). Denote also $f(t)=\frac{1}{t}\alpha(t)$,
\[g(t)=t^{\rho -1} h_n(1/t,\theta_0,q),\]
and the integral in (2.2) as $J(r)$. Integrating $J(r)$ by parts, we get
\[J(r)=-\frac{1}{r}\int_0^{\infty} f(t) \; t \frac{\partial}{\partial t}g\left(\frac{t}{r}\right) dt, \]
since the integrated term vanishes due to (2.1).
We have
\[f(t)=\frac{\alpha(t)}{t}=\frac{1}{t}\int_0^t s^{3-n-\rho}d(s^{n-2}n(t)).\]
For a subharmonic function $u$ of order $\rho$ and finite type, its counting function $n(t)$ is non-decreasing and satisfies $n(t)\leq Ct^{\rho}$. Therefore, after integrating $f(t)$ by parts, we get
\beq f(t)=\frac{\alpha(t)}{t}=\frac{n(t)}{t^{\rho}}+\frac{n-3+\rho}{t}\int_0^t\frac{n(t)}{t^{\rho}} dt. \eeq
Since
\[\frac{n(y)}{y^{\rho}}-\frac{n(x)}{x^{\rho}}=\frac{n(y)-n(x)}{x^{\rho}}+\frac{n(y)}{y^{\rho}} \left(1-\left(\frac{y}{x}\right)^{\rho}\right), \]
the integrated term in (2.3) is a slowly decaying function. The same is valid for the integral in (2.3), since
\[\frac{1}{y}\int_0^y \frac{n(s)}{s^{\rho}}ds - \frac{1}{x}\int_0^x \frac{n(s)}{s^{\rho}}ds=\frac{1}{y}\int_x^y \frac{n(s)}{s^{\rho}}ds +\left(\frac{x}{y}-1\right) \frac{1}{x}\int_0^x \frac{n(s)}{s^{\rho}}ds \]
\[\leq C\frac{y-x}{y} + C\frac{|x-y|}{y} \rightarrow 0 \]
as $x\rightarrow \infty$ and $y/x \rightarrow 1$.   \\

To apply Theorem W, we write the integral $J(r)$ as
\[J(r)= \frac{1}{r}\int_0^{\infty} \left(-u \frac{\partial}{\partial u}g(u) \right)\biggl|_{u=t/r} f(t) dt, \]
and we have to know the location of the zeros of the Mellin transformation of the kernel in (2.2), that is, of the integral
\beq M(g;v)=-\int_0^{\infty} u\frac{\partial}{\partial u} g(u) u^{-\imath v}du=-\int_0^{\infty} \frac{\partial}{\partial u} g(u) u^{1-\imath v}du \eeq
as a function of $v$. Integrating (2.4)  by parts, we again notice that the integrated term vanishes, and we have
\[M(g;v)=(1-\imath v) \int_0^{\infty} g(u) u^{-\imath v}du. \]

We remind that $g(u)=u^{\rho-1}h_n(1/u,\phi,q)$, thus we consider the integral
\[M(h_n;v)=(1-\imath v)\int_0^{\infty} u^{\rho-1}h_n\left(\frac{1}{u}, \phi,q\right) u^{-\imath v}du \]
\[=(1-\imath v)\int_0^{\infty} h_n(u, \theta_0,q) u^{-1-\rho-\imath v}du, \]
and by virtue of (A.2)-(A.4) with $s=-\rho-\imath v$, we express $M(h_n;v)$ through the associated Legendre functions of the first kind
\[P^{\frac{3-n}{2}}_{-\rho -\imath v -\frac{n-1}{2}}(\cos \phi).\]

As we have stated above, the latter has no complex roots; while if $v=0$, it has $E(\rho+1)$ zeros, which are excluded by the condition $\phi \not\in \Theta_n(\rho)$, thus
\[M(h_n;0)=\frac{\pi 2^{\frac{n-3}{2}} \Gamma(\frac{n-1}{2}) \prod_{k=1}^{n-2}(\rho+k) \mathbf{\Delta}}{(n-3)! \sin (\pi \rho )
\left(\sin \phi\right)^{\frac{n-3}{2}}}P^{\frac{3-n}{2}}_{-\rho-\frac{n-1}{2}}(\cos \phi). \]

After some simple algebra we derive equation (1.12).  \\

To show that the restriction $\phi \not\in\Theta_n(\rho)$ is essential, we consider, in the simplest case $n=3$ and $0<\rho<1$, the function
\[u_0(x)=r^{\rho}+r^{\rho} \sin \log \log r \; P_{\rho}(\cos \theta_1),\]
where as usual, $P_{\rho}=P^0_{\rho}$. We straightforwardly verify that the Laplacian
\[\Delta u_0(x)\]
\[=r^{\rho -2}\left\{\rho (\rho+1)+(2\rho+1)\frac{\cos \log \log r}{\log r}P_{\rho}(\cos \theta_1)+ O\left(\frac{1}{\log ^2 r}\right) \right\}.\]
It is known \cite[Sect.3.9.2 (15)]{BE} that $P_{\rho}(x)\approx \frac{\sin \pi \rho}{\pi} \log \frac{1+x}{2}$ as $x\rightarrow -1$, thus $\Delta u_0(x)\geq 0$ for $r>r_0$ uniformly in $\theta_1$. Whence, the function $u_0$ is subharmonic in $\mathbb{R}^3$ outside of a fixed ball $B(r_0)$ if $r_0$ is large enough. Moreover, if $\phi \not\in \Theta_3(\rho)$, then limits (1.11) do not exist since $u_0$ is oscillating. On the other hand, for $\phi \in \Theta_3(\rho)$, limit (1.11) does exist since $P_{\rho}(\cos \phi)=0$.        $ \hfill \qed$  \\

\noindent\textbf{Proof of Theorem 2}. As before, we assume $u(0)=0$ and represent the function $u$ as \cite[Eq-n 4.5.15]{HaKe},
\beq u(x)=\int_0^{\infty} \frac{P_n(|x|,t,\theta_1) N(t) dt}{\left(|x|^2+2|x|t\cos \theta_1 +t^2 \right)^{n/2+1}}. \eeq
In this proof we use the following special case of the Tauberian theorem of Drasin \cite{Dra1}:

\emph{Let} $f(t),\; -\infty <t< \infty$, \emph{be a positive increasing function such that its convolution with a kernel} $k(t)$
satisfies
\[\int_{-\infty}^{\infty}k(t-y)f(y)dy=\left\{L+o(1)\right\} f(t),\; \tti .\]
\emph{If the kernel} $k$ \emph{is positive almost everywhere, and its Laplace transform }
\[\mathcal{L}k(s)=\int_{-\infty}^{\infty} e^{-st} k(t) dt \]
\emph{exists if} $-\sigma <s<\varrho$ \emph{for some positive, maybe infinite} $\sigma$ \emph{and} $\varrho$, \emph{whereas}
$\mathcal{L}k(-\sigma)=\mathcal{L}k(\varrho)=\infty$, then
\[f(t)=e^{\lambda t} \psi(t), \]
\emph{where for every fixed} $a$, $\psi(t+a)/\psi(t)\rightarrow 1, \; \tti$. \emph{Moreover}, $\lambda (\geq 0)$ \emph{must satisfy} $\mathcal{L}k(\lambda)=L$.   \\

After substituting $t=e^y$ and $r=|x|=e^p$, (2.5) becomes
\[u(e^p,\theta_1)=\int_{-\infty}^{\infty}k(p-y)f(y)dy,\]
where $f(y)=N(e^y)$ and
\[k(t)=\frac{e^t\left\{(n-1)\cos \theta_1 +[n+(n-2)\cos^2 \theta_1]e^t+(n-1)\cos \theta_1 e^{2t}\right\}}{\left(1+2 e^t \cos \theta_1 +e^{2t}\right)^{n/2+1}}.\]
It is clear that $k(t)\geq 0$ for all $0\leq \theta_1 \leq \pi/2$. Moreover, elementary considerations with $P_n(r,t,\theta_1)$ as a quadratic trinomial with respect to $\cos \theta_1$, show that the right boundary $\theta_1 \leq \pi/2$ cannot be increased. In particular, $\sigma =-n-1$, $\varrho=1$, and all the other conditions also easily verified.

We write down the final result only in the case $\theta_1 =0$. The kernel is
\[k(t)=(n-1)e^{(n-1)t}\left(1+e^t\right)^{-n},\]
thus its positivity is obvious, the computation of its Laplace transform is straightforward, and the transcendental equation for the order $\rho$ is
\[\overline{\mathbf{\Delta}}=\frac{ \Gamma(n-1-\rho)}{(n-2)! \Gamma(1-\rho)} \frac{\pi \rho}{\sin (\pi \rho)}.\]
Returning to the variables $r=e^p$ and $t=e^y$, we get for this $\rho$,
\[N(t)=t^{\rho}\psi_1(t)\]
and
\[u(r)=r^{\rho} \psi(r),\]
where $\psi$ and $\psi_1$ are slowly varying functions. Finally, similarly to Proposition 1, we can straightforwardly derive the asymptotic representation of $u(r,\theta)$ for any $\theta$.  $\hfill \qed$

\appendix
\section{Mellin Transform of the Riesz Kernel}

The associated Legendre functions of the first kind $P^{\mu}_{\nu}(z)$ are particular solutions of the Legendre differential equation \cite[Chap. 3]{BE}
\[(1-z^2)y''(z)-2zy'+\left[\nu (\nu +1)-\mu^2 (1-z^2)^{-1}\right] y=0,\]
where $\mu$ and $\nu$ are, in general, complex parameters; $P^{\mu}_{\nu}(\xi)$ stands for these functions on the cut $-1<\xi <1$ \cite[Sect. 3.4]{BE}. In many instances the Legendre functions are a natural replacement of the trigonometric functions in many-dimensional problems, so that there is a continuous stream of research regarding the spherical functions. In particular, it may be of interest to represent certain analytic objects, like series and integrals as Legendre's functions; among the newest publications see, e.g., \cite{Szm}. The following statement is used in the proof of the main result of this paper.

Instead of the Weierstrass canonical kernel in (1.4), we consider a slightly more general case of the Riesz kernel
\[k_{\lambda}(t,\xi)=(1+t^2+2t\xi)^{-\lambda} \]
with any $\lambda$ such that $\Re \lambda >0$. It is known that its Mellin transform can be represented through the associated Legendre's function,
\beq \begin{array}{cc} \int_0^{\infty} t^{\nu-\mu}\left(1+t^2+2t\xi \right)^{\mu-1/2}dt \vspace{.3cm} \\
=\frac{\Gamma(1-\mu) \Gamma(\nu-\mu+1) \Gamma(-\mu-\nu)}{2^{\mu} \Gamma(1-2\mu)} \left(1-\xi^2\right)^{\mu/2} P^{\mu}_{\nu}(\xi) , \end{array}   \eeq
where $\Gamma$ is Euler's $\Gamma-$function, \cite[Sect. 6.2, Eq. (22)]{BaEIT}. The integral in (A.1) is convergent if
\[ \Re \mu - \Re \nu <1 \mbox{ and } \Re \mu +\Re \nu <0. \]
It should be also mentioned that the kernel $k_{\lambda}$ is a generating function of the Gegenbauer polynomials $G^{\lambda}_j (\xi)$, Cf. the case $\lambda=(n-2)/2$ in Section 1,
\[k_{\lambda}(t, \xi)=(1+2t\xi+ t^2)^{-\lambda} = \sum_{j=0}^{\infty} (-t)^j G^{\lambda}_j (\xi). \]

We consider only the case $n\geq 3$, since for $n=2$ the spherical functions are essentially the trigonometric functions, and the integral (A.1) for $n=2$ is well known, it is equal to $\frac{\pi \cos \rho \theta}{\rho \sin \pi \rho}$. Thus, let
\[h(u)=h(\lambda, q, u, \xi) = -\left(1+2u \xi +u^2\right)^{-\lambda} + \sum^q_{j=0} (-u)^j G_j^{\lambda}(\xi)  \]
\[=-k_{\lambda}(u, \xi)+ \sum^q_{j=0} (-u)^j G_j^{\lambda}(\xi),\]
which has the same bound as (2.1),
\[|h(\lambda, q, u, \xi)|\leq C \min \{u^q; u^{q+1}\},\; 0<u<\infty, \]
where a positive constant $C$ does not depend on $u$ and $\xi$.   \\

\noindent\textbf{Proposition}. \emph{For any integer} $q=0,1,2,\ldots $, \emph{real} $\xi$, $-1< \xi <1$, \emph{and complex} $\lambda$ \emph{and} $s$ \emph{such that}
\[0< \Re \lambda \]
\emph{and}
\[-q-1< \Re s <-q, \]
\emph{the Mellin transformation of} $h \vspace{.2cm}$ \emph{is}
\beq \begin{array}{cc} M(h,s)=\int_0^{\infty} \left\{-\left(1+u^2+2u\xi\right)^{-\lambda} + \vspace{.8cm} \sum_{j=0}^q (-u)^j G_j^{\lambda}(\xi)\right\}u^{s-1}du  \vspace{.4cm} \\
=-\frac{\sqrt{\pi}\Gamma(s) \Gamma(2\lambda-s)}{2^{\lambda-1/2}\Gamma(\lambda)} \left(1-\xi^2\right)^{\frac{1-2\lambda}{4}}
P_{s-\lambda-1/2}^{1/2-\lambda}(\xi). \end{array} \eeq

\noindent\textbf{Proof}. The Mellin transform of the kernel $h$,
\beq M(h,s)=\int_0^{\infty} h(u)u^{s-1}du \eeq
is convergent for $-1-q<\Re s<-q$. We integrate by parts the integral in (A.3) $q+1$ times, so that the polynomial part of $h$ vanishes, whence
\[\frac{\partial^{q+1}}{\partial u^{q+1}}k_{\lambda}(u, \xi)=-\frac{\partial^{q+1}}{\partial u^{q+1}}h(u)\]
and
\beq M(h,s)=\frac{(-1)^q}{\prod_{k=0}^q(s+k)}\int_0^{\infty}u^{s+q}\frac{\partial^{q+1}}{\partial u^{q+1}} \left((1+u^2+2u\xi)^{-\lambda} \right) du. \eeq
The latter integral is convergent for non-integer $s$ such that
\[-1-q<\Re s<2\Re \lambda,  \]
thus providing the analytic continuation of M(h,s) as a meromorphic function into this broader domain of the $s-$plane.  \\

Now we consider the Mellin transform of the kernel $k_{\lambda}$,
\[ M(k_{\lambda},s)=\int_0^{\infty}(1+t^2+2t\xi)^{-\lambda}t^{s-1} dt, \]
which is convergent for
\beq  -q-1<\Re s< 2\Re \lambda. \eeq
Integrating it by parts $q+1$ times, we get
\beq M(k_{\lambda},s)=\frac{(-1)^{q+1}}{\prod_{k=0}^q(s+k)}\int_0^{\infty}t^{s+q}\frac{\partial^{q+1}}{\partial t^{q+1}} \left((1+t^2+2t\xi)^{-\lambda}\right) dt. \eeq
Due to (A.5), all the integrated terms vanish and the integral is convergent in the wider region $-1-q<\Re s<2\Re \lambda$ excluding the poles at integer points.   \\

By (A.1), we express $M(k_{\lambda},s)$ through $P^{\mu}_{\nu}(\xi)$ with $\mu=\frac{1}{2}-\lambda$ and $\nu=s-\frac{1}{2}-\lambda$. Finally, combining (A.4) and (A.6) and using Legendre's formula for the $\Gamma-$function with double argument,
\[\Gamma(2z)=2^{2z-1} \pi^{-1/2}\Gamma(z)\Gamma(1/2+z),\]
we get the result.

The Legendre functions $P^{\mu}_{\nu}(\xi)$ are entire functions in both $\mu$ and $\nu$, thus all the equations are justified due to the principle of analytic continuation.  \hfill  $\qed$  \\

\noindent\textbf{Remark}. If $q=0$, a simpler proof can be given. In this case, we can explicitly compute the derivative in (A.4) and split the integral into the two integrals of kind
\[\int_0^{\infty}u^{\alpha}(1+u^2+2u\xi)^{\beta} du  \]
with two different $\alpha$. Then we apply (A.1) to express each of them as $P^{\mu}_{\nu}(z)$ and use the recurrence formula \cite[Sect. 3.8, Eq-n (18)]{BE}
\[\begin{array}{rr} (\nu-\mu+1)P_{\nu+1}^{\mu}(\cos \theta_1)-(\nu+\mu+1)\cos \theta_1 P_{\nu}^{\mu}(\cos \theta_1) \vspace{.3cm}\\
          =\sin \theta_1 P_{\nu}^{\mu+1}(\cos \theta_1) \end{array} \]
to arrive at (A.2) with $q=0$. However, the explicit computation of the derivatives for bigger $q$ becomes cumbersome.  \\

In the case, we are interested in, $\lambda=\frac{n-2}{2}$ and $s=-\rho$, and formula (A.2) reads as follows. \\

\noindent\textbf{Corollary}. For any integer $q=0,1,2,\ldots $, real $\xi$, $-1< \xi <1$, and complex $\rho$ such that
\[q< \Re \rho < q+1, \]
the Mellin transform of the function $h_n$ is
\[M(h_n,\rho)=-\frac{\pi \sqrt{\pi} 2^{(3-n)/2} \prod_{k=1}^{n-3}(\rho+k)(1-\xi^2)^{(3-n)/4}}{\sin \pi \rho \; \Gamma((n-2)/2)}
P_{-\rho-(n-1)/2}^{(3-n)/2}(\xi), \]
where $\xi =\cos \theta_1$, or using the equation
\[\sqrt{\pi} (n-3)! =2^{n-3}\Gamma((n-1)/2)\Gamma((n-2)/2),\]
which can be immediately proved by induction,
\[M(h_n,\rho) \vspace{.2cm}  \]
\[=\frac{\pi 2^{(n-3)/2} \prod_{k=1}^{n-3}(\rho+k)\Gamma((n-1)/2)(1-\xi^2)^{(3-n)/4}}{(n-3)! \sin \pi \rho }
P_{-\rho-(n-1)/2}^{(3-n)/2}(\xi). \]

\bigskip

\textbf{Acknowledgement} The author is grateful to Professors P.Z. Agranovich, T.M. Dunster, and F.W.J. Olver for very useful comments.

This work was done during a sabbatical leave from the City University of New York.

\smallskip

\end{document}